\documentclass[12pt]{amsart}
\usepackage{txfonts}
\usepackage{subfig}
\usepackage{amscd,amssymb,mathabx}
\usepackage[arrow,matrix,graph,frame,poly,arc,tips]{xy}
\usepackage{graphicx,calrsfs}
\usepackage{mathtools}
\usepackage{tikz}
\usetikzlibrary{decorations.markings}
\usetikzlibrary{shapes}
\usetikzlibrary{backgrounds}
\usetikzlibrary{calc}
\usepackage{mdwlist}
\usepackage{multirow}
\usepackage{enumerate}
\usepackage{verbatim}

\DeclarePairedDelimiter{\form}{\langle}{\rangle}

\newcommand\ba{\begin{align*}}
\newcommand\ea{\end{align*}}
\newcommand\be{\begin{enumerate}}
\newcommand\ee{\end{enumerate}}
\newcommand\bp{\begin{proof}}
\newcommand\ep{\end{proof}}
\newcommand\bpp{\begin{prop}}
\newcommand\epp{\end{prop}}
\newcommand\bpb{\begin{prob}}
\newcommand\epb{\end{prob}}
\newcommand\bd{\begin{defn}}
\newcommand\ed{\end{defn}}
\newcommand\bh{\begin{hint}}
\newcommand\eh{\end{hint}}

\newcommand\bZ{\mathbb{Z}}

\newcommand\supp{\operatorname{supp}}
\newcommand\Id{\operatorname{Id}}
\newcommand\lk{\operatorname{Lk}}
\newcommand\st{\operatorname{St}}

\newcommand\sse{\subseteq}
\newcommand\co{\colon}

\def\thetitle{{Finite subgraphs of an extension graph}}
\def\theauthors{{Sang-hyun Kim and Thomas Koberda}}
\usepackage{hyperref}
\hypersetup{
   colorlinks=false,
   plainpages,
   urlcolor=black,
   linkcolor=black
   pdftitle=   \thetitle,
   pdfauthor=  {\theauthors}
}

\theoremstyle{theorem}
\newtheorem{thm}{Theorem}[section]
\newtheorem{lem}[thm]{Lemma}

\newtheorem{prop}[thm]{Proposition}

\newtheorem*{claim*}{Claim}

\theoremstyle{remark}
\newtheorem{exmp}[thm]{Example}

\theoremstyle{definition}
\newtheorem{defn}[thm]{Definition}
\newtheorem{prob}{Problem}[section]

\begin{document}
\title\thetitle
\date{\today}

\author[S. Kim]{Sang-hyun Kim}
\address{Department of Mathematical Sciences, Seoul National University, Seoul, Korea}
\email{s.kim@snu.ac.kr}
\urladdr{http://cayley.kr}

\author[T. Koberda]{Thomas Koberda}
\address{Department of Mathematics, University of Virginia, Charlottesville, VA 22904-4137, USA}
\email{thomas.koberda@gmail.com}
\urladdr{http://faculty.virginia.edu/Koberda}

\author[J. Lee]{Juyoung Lee}
\address{Department of Mathematical Sciences, Seoul National University, Seoul, Korea}
\email{ljy219@snu.ac.kr}

\begin{abstract}
Let $\Gamma$ be a finite graph and let $\Gamma^{\mathrm{e}}$ be its extension graph.
We inductively define a sequence  $\{\Gamma_i\}$ of finite induced subgraphs of $\Gamma^{\mathrm{e}}$ through successive applications of an operation called ``doubling along a star''.
Then we show that every finite induced subgraph of $\Gamma^{\mathrm{e}}$ is isomorphic to an induced subgraph of  some $\Gamma_i$.
This result strengthens \cite[Lemma 3.1]{KK2013GT}.
\end{abstract}

\maketitle

\section{Statement of the result}
\subsection{Notations}
\emph{Throughout this note}, let us fix a finite graph $\Gamma$
and its vertex set
\[
V = V(\Gamma) = \{a_0,\ldots,a_{n-1}\}.\]
We will mostly follow the terminology and notations in~\cite{KK2013GT}.

For $U,W\sse A(\Gamma)$, we define
\[ U^W =\{u^w\co u\in U\text{ and }w\in W\}\sse A(\Gamma).\]
We put 
  \[
V^{\mathrm{e}} = V^{A(\Gamma)}.\] 
Recall that the \emph{extension graph} $\Gamma^{\mathrm{e}}$ is defined as the commutation graph of $V^{\mathrm{e}}$ in $A(\Gamma)$; see Definition~\ref{defn:universal}.

It will be convenient for us to denote
\[ V^{\mathrm{e}}_\bZ= \{(v^k)^g\co v\in V, \  k\in\bZ\setminus0,\ g\in A(\Gamma)\}\sse A(\Gamma).\]
We have a map $(\cdot)^\ast\co V^{\mathrm{e}}_\bZ\to V^{\mathrm{e}}$ 
 defined by the formula
\[
\left((v^k)^g\right)^* = v^g.\]
For example, if $a,b,c\in V$ then 
we have $(a^{-2})^{bc})^*=a^{bc}$.

For each $w\in A(\Gamma)$, we let $\|w\|$ denote the \emph{word length} of $w$. 
In other words, $\|w\|$ is the smallest nonnegative integer $\ell$ such that we can write
\[ w = s_1^{e_1}\cdots s_\ell^{e_\ell}\]
for some  $s_i\in V$ and $e_i\in\{-1,1\}$. 
In this case, we define the \emph{support} of $w$ as 
\[\supp w = \{s_1,s_2,\ldots,s_\ell\}\sse V.\]


\subsection{Double of a graph along a star}
Let $X$ be a graph. For $S\sse V(X)$, we denote by $X(S)$ or by $XS$ the subgraph of $X$ induced by $S$.
The \emph{star} of a vertex $v$ in $X$ is the set of vertices in $X$ that are either equal or adjacent to $v$. We denote the star of $v$ as $\st_X(v)$ or $\st(v)$. We define the \emph{link} of $v$ as
\[
\lk_X(v) = \st_X(v)\setminus\{v\}.\]

Fix a vertex $v$ of $X$. 
Let $S_0,S_1$ be sets with some fixed bijections \[ \rho_i\co S_i\to  V(X)\setminus\st(v).\]
Then we can define a new graph $Y$ by requiring that
\[ V(Y) = S_0\coprod S_1\coprod\st (v)\] and that $\{a,b\}\in E(Y)$ if and only if  one of the following holds:
\be[(i)]
\item
$a,b\in S_i$ and $\{\rho_i(a),\rho_i(b)\}\in E(X)$ for some $i=0$ or $1$.
\item
$a\in S_i,b\in \st(v)$ and $\{\rho_i(a),b\}\in E(X)$ for some $i=0$ or $1$.
\item
$a,b\in \st(v)$ and $\{a,b\}\in E(X)$.
\ee
The graph $Y$ thus obtained is called the \emph{double of $X$ along the star of $v$.}

The main result of this paper is the following.

\begin{thm}[compare with {\cite[Lemma 3.1]{KK2013GT}}]\label{t:double}
There exists an infinite sequence of finite induced subgraphs $\{\Gamma_i\}_{i\ge0}$  of $\Gamma^e$
such that 
$\Gamma_{i+1}$ is the double of $\Gamma_i$ along a star
and such that
every finite induced subgraph of $\Gamma^{\mathrm{e}}$
admits an embedding into some $\Gamma_i$ as an induced subgraph.
 \end{thm}
 
This theorem strengthens Lemma 3.1 in~\cite{KK2013GT}, where the proof of the lemma is omitted. 
As there has been much interest recently concerning on the combinatorial structures of extension graphs~\cite{CDK2013JA,CasalsRuiz2015IMRN,LL2017arXiv,Huang2014qi,Huang2016qi2,HK2016},
we decided to write down a very detailed construction of such a sequence $\{\Gamma_i\}_{i\ge0}$.

\section{Universal sequence}
We denote the symmetric difference of two sets $A$ and $B$ as $A\triangle B$.
\bd\label{defn:universal}
Let $X$ be a subset of a group $G$. The \emph{commutation graph} of $X$, denoted as $\operatorname{CG}(X)$, is the simplicial graph 
whose vertex set is $X$
and in which two distinct vertices $x,y\in X$ are adjacent if and only if $[x,y]=1$.
\ed
Note that we have the \emph{natural homomorphism}
\[
A(\operatorname{CG}(X))\to \form{X}\]
defined by the unique extension of $\Id_X$.

\begin{lem}[\cite{BKS2008,KK2013GT}, \emph{cf. }\cite{Kim2008,Bell2011}]\label{l:double-raag}
Let $G$ be a group and $X_0\sse G$ be a subset such that the natural homomorphism
\[
A(\operatorname{CG}(X_0))\to \form{X_0}\]
is an isomorphism.
Suppose $u\in X_0$, and define
\[ X_1 = \left( X_0\cup X_0^u\right)\triangle \{u,u^{-2}\}.\]
We define $\phi\co \form{X_0}\to\bZ_2$ by 
 $\phi(u)=1$ and $\phi(v)=0$ for all $v\in X_0\setminus\{u\}$.
Then $\operatorname{CG}(X_1)$ is the double of 
 $\operatorname{CG}(X_0)$ along the star of $u\in X_0$.
Moreover, we have the commutative diagram
\[
\xymatrix{
1 \ar[r] & A(\operatorname{CG}(X_1))\ar[r]\ar[d]^\cong&
A(\operatorname{CG}(X_0))\ar[r]\ar[d]^\cong& \bZ_2\ar[r]\ar[d]^\Id & 1 \\
1 \ar[r] & \form{X_1}\ar[r]&
\form{X_0}\ar[r]^{\phi}& \bZ_2\ar[r] & 1 \\}\]
where the left two vertical maps are natural isomorphisms.
\end{lem}

Let us now define an infinite sequence $\{u_i\}_{i\ge0}\sse V^{\mathrm{e}}_\bZ$ as follows:
\begin{align*}
&u_0=a_0,\ 
u_1 = a_0^{-2},\
u_2 = a_1,\
u_3 = a_1^{-2},\
\ldots,\ 
u_{2n-1}=a_{n-1}^{-2},\\
&u_{2n} = a_0^4,\ u_{2n+1} = a_0^{-8},\ u_{2n+2}=a_1^4,\ \ldots,\ 
u_{4n-1}=a_{n-1}^{-8},\\
&u_{4n} = a_0^{16},\  u_{4n+1} = a_0^{-32},\ \ldots\text{ and so forth.}
\end{align*}

For each $i\ge0$ and $0\le j< 2n$,
one can more succinctly write
\[u_{2ni+j} = a_{\lfloor{j/2}\rfloor}^{r_{i,j}},\quad r_{i,j}={(-2)^{2i+j-2\lfloor j/2\rfloor}}.\]

We then have an infinite sequence $\{U_i\}_{i\ge0}$ of subsets of $V^{\mathrm{e}}_\bZ$ defined as
\begin{align*}
U_0&=V,\\
U_{i+1}&=\left(U_{i}\cup U_{i}^{u_i}\right)\triangle\{u_i,u_i^{-2}\}.
\end{align*}
We will call the sequence $\{(u_i,U_i)\}_{i\ge0}$ as a \emph{universal sequence} in $\Gamma^e$. Note that a universal sequence depends on the choice of the enumeration $V=\{a_0,\ldots,a_{n-1}\}$.

The following lemma is the basic building block for our construction.
\begin{lem}\label{l:univ}
The following hold for each $i\ge0$.
\be
\item $u_i\in U_i$.
\item
The natural homomorphism $A(\operatorname{CG}(U_{i}))\to\form{U_i}$ is an isomorphism,
and moreover,
$\operatorname{CG}(U_{i+1})$ is the double of 
$\operatorname{CG}(U_{i})$ along the star of $u_i$.
\item
The map $x\mapsto x^*$ is injective on $U_i$.
\item
The two graphs $\operatorname{CG}(U_{i})$ and $\Gamma^e U_i^*$ are  isomorphic by the isomorphism $x\mapsto x^*$.
\ee
\end{lem}

\bp For $i\ge0$ and $0\le j<n$ let us note
\begin{align*}
u_{2ni+2j}= a_j^{4^i},\quad & \{a_0,\ldots,a_{j-1}\}^{4^{i+1}}\cup
\{a_j,\ldots,a_{n-1}\}^{4^i}\sse U_{2ni+2j}.\\
u_{2ni+2j+1}= a_j^{-2\cdot 4^i},\ & 
\{a_0,\ldots,a_{j-1}\}^{4^{i+1}}\cup
\{a_j\}^{-2\cdot 4^i}\cup
\{a_{j+1},\ldots,a_{n-1}\}^{4^i}\sse U_{2ni+2j+1}.
\end{align*}
So, part (1) is obvious. 
Part (2) follows from an induction combined with Lemma~\ref{l:double-raag}.

For part (3), 
assume $x,y\in U_i$ satisfy $x^*=y^*$.
Then we can write $x=(u^g)^k$ and $y=(u^h)^m$ for some $u\in V$, $g,h\in A(\Gamma)$ and $k,m\in\bZ\setminus0$.
Since $x^m = y^k$ in $A(\Gamma)$,  part (2) implies that $x=y$.

Consider  part (4). The map $(\cdot)^*$ defines a natural bijection 
\[
\operatorname{CG}(U_i)\to \Gamma^e U_i^*.\]
Note that for $x,y\in V^e_\bZ$, we have
\[ [x,y]=1\Leftrightarrow [x^*,y^*]=1.\]
Hence $\operatorname{CG}(U_{i})$
and
$\Gamma^eU_i^*=\operatorname{CG}(U_{i}^*)$ are isomorphic.
\ep

The following lemma shows that the sequence $\{\operatorname{CG}(U_i)\}_{i\ge0}$ eventually contains copies of all the finite induced subgraphs of $\Gamma^{\mathrm{e}}$. 

\begin{lem}\label{l:sigma}
For each finite set $W\sse A(\Gamma)$, there exists $K>0$ and a map
\[
\sigma\co W\to A(\Gamma)\]
such that the following hold:
\be[(i)]
\item
$\Gamma^e(V^W) \cong \Gamma^e(V^{\sigma W})$.
\item
$V^{\sigma W}\sse U_K^*$.
\ee
\end{lem}

The proof of this lemma is postponed until the next section.
Let us first deduce the main theorem of this paper.

\bp[Proof of Theorem~\ref{t:double} assuming Lemma~\ref{l:sigma}]
We can find a finite set $W\sse A(\Gamma)$ such that $V(\Lambda)\sse V^W$.
By the conditions (i) and (ii) of the lemma, the graph $\Lambda$ is an induced subgraph of 
$\Gamma^e U_K^*$. Lemma~\ref{l:univ} implies that 
 $\Gamma^e U_K^*$ is obtained from $\Gamma$ by successive applications of doubling along stars, as desired.
\ep

\section{Inflating powers of letters}
In this section, we find a map $\sigma$ satisfying the conditions of Lemma~\ref{l:sigma}.
\subsection{Canonical expression}
Consider an arbitrary $w\in A(\Gamma)$. We can write
\begin{equation}
w= s_1^{e_1} s_2^{e_2}\cdots s_\ell^{e_\ell}\tag{*}\label{eqn:canonical}
\end{equation}
where $\ell=\|w\|$,  $s_i\in V$ and  $e_i\in\{-1,1\}$.
For each $i=1,2,\ldots,\ell$, we put
\[f_i = \min_y \|y\|\]
where $y$ varies among the words in $A(\Gamma)$ such that we can write
 \[s_i^{e_i} s_{i+1}^{e_{i+1}}\cdots s_\ell^{e_\ell} = x \cdot s_i^{e_i}\cdot y\]
for some word 
$x\in\form{\lk_\Gamma(s_i)}$.
Roughly speaking, we minimize the length of a word $y$ that ``remains'' on the right of $s_i^{e_i}$.
We call $(f_1,\ldots,f_\ell)$ as the \emph{right-counting vector} corresponding to the word in \eqref{eqn:canonical}.

\begin{exmp}\label{exmp:canonical}
Let $V=\{a_0=a,\ a_1=b,\ a_2=c\}$ and 
\[A(\Gamma)=\form{a,b,c\mid [a,b]=1},\quad w=a^2b^{-3}cb.\]
The right-counting vector for this word is 
 \[(f_1,\ldots,f_7)=(3,2,4,3,2,1,0).\]
Let us consider a different word representing the same element:
\[
w = b^{-1}ab^{-1}ab^{-1}cb.\]
Then the corresponding right-counting vector becomes
\[(4,3,3,2,2,1,0).\]
\end{exmp}


A word in~\eqref{eqn:canonical} is called a \emph{canonical expression for $w$} if the corresponding right-counting vector satisfies the following two conditions:
\be[(A)]
\item
$f_1\ge f_2\ge\cdots\ge f_\ell$;
\item
If $f_i = f_j$ for some $i<j$ and if we write $s_i = a_p$, $s_j = a_q$ for some $0\le p,q<n$,
then we have that $p<q$ and that $[a_p,a_q]=1$.
\ee

\begin{lem}\label{l:canonical-expr}
Each $w\in A(\Gamma)$ admits a unique canonical expression.
\end{lem}
\bp
Let $(f_1,\ldots)$ be the right-counting vector for a word representing $w$. Put
\begin{align*}
A&=(\#\text{ of }(i,j)\text{ where }i<j\text{ and }f_i<f_j),\\
B&=(\#\text{ of }(i,j)\text{ where }i<j,\ f_i=f_j\text{ and }s_i=a_p,s_j=a_q\text{ for some }p>q).
\end{align*}
A canonical expression for $w$ is then obtained by minimizing the lexicographical order of the tuple $(A,B)$, resulting in $(0,0)$.
\ep

\subsection{Proof of Lemma~\ref{l:sigma}}
Recall our notation $\{(u_i,U_i)\}_{i\ge0}$, which is defined in the previous section.
In order to prove Lemma~\ref{l:sigma}, it suffices for us to consider the case when
\[W =B(M)= \{w\in A(\Gamma)\co \|w\|\le M\}\] for some positive integer $M$. 

Let $w\in B(M)$ be written as~\eqref{eqn:canonical}, which is not necessarily canonical.
Denote by $(f_1,\ldots,f_\ell)$ the corresponding right-counting vector. 
We define 
\begin{align*}
&N_i =  \frac{3e_i-1}2 \cdot 4^{M-1-f_i}\in \{\pm 2^m\co m\ge0\},\\
&\sigma( s_1^{e_1} s_2^{e_2}\cdots s_\ell^{e_\ell})= s_1^{N_1} s_2^{N_2}\cdots s_\ell^{N_\ell}.\tag{**}\label{eqn:sigma}
\end{align*}

\begin{lem}\label{l:sigma-defn}
The following hold.
\be
\item
The map $\sigma\co B(M)\to A(\Gamma)$ is well-defined. That is, if two words $x_1$ and $x_2$ represent the same element $w$ in $A(\Gamma)$, then $\sigma(x_1)=\sigma(x_2)$ in $A(\Gamma)$.
\item
If \eqref{eqn:canonical} is a canonical expression for $w\in B(M)$
and if $\sigma(w)$ is written as \eqref{eqn:sigma},
then $\left\{s_i^{N_i}\right\}_{1\le i\le\ell}$ is a subsequence of $\{u_i\}_{i\ge0}$.\ee
\end{lem}

\bp
(1)
One goes through the definition of $\sigma$ for a different expression
\[
w = s_1^{e_1}\cdots s_{i-1}^{e_{i-1}} s_{i+1}^{e_{i+1}} s_{i}^{e_{i}} s_{i+2}^{e_{i+2}}\cdots s_\ell^{e_\ell}\]
when $[s_i,s_{i+1}]=1$ and verifies that the resulting element coincides with $\sigma(w)$.

(2)
Immediate from the definition of a right-counting vector and a canonical expression.
\ep

\begin{exmp}
Continuing Example~\ref{exmp:canonical} and setting $M=7$, we have
\[
\sigma(w)=\sigma(b^{-1}ab^{-1}ab^{-1}cb) = b^{-2\cdot 4^2} a^{4^3} b^{-2\cdot 4^3} a^{4^4} b^{-2\cdot 4^4} c^{4^5}b^{4^6}.\]
\end{exmp}

We claim the map $\sigma\co B(M)\to A(\Gamma)$ thus defined satisfies the conditions of Lemma~\ref{l:sigma}. 
The condition (ii) is implied by Lemma~\ref{l:sigma-defn} (2), so it remains to show the condition (i).
For $v\in V$, we let $Z(v)=\form{\st_\Gamma(v)}$, which is the centralizer group of $v$ in $A(\Gamma)$.
\begin{lem}\label{l:uv}
For $u,v\in V$ and $x,y\in B(M)$, we have the following.
\be
\item  $xy^{-1}\in Z(u)$ iff $\sigma(x)\sigma(y)^{-1}\in Z(u)$.
\item  $xy^{-1}\in Z(v) Z(u)$ iff $\sigma(x)\sigma(y)^{-1}\in  Z(v)Z(u)$.
\ee\end{lem}

\bp
Let us consider reduced expressions
\[x = x_0\cdot  p,\quad y = y_0\cdot  p\]
such that $x_0\cdot y_0^{-1}$ is reduced.
Then we have reduced expressions
\[
\sigma(x) = x_1\cdot \sigma(p),\quad\sigma(y)=y_1\cdot \sigma(p)\]
for some words $x_1$ and $y_1$.
We can write
\begin{align*}
x_0 y_0^{-1} &= \prod_{i=1}^k t_i^{g_i},\\
x_1 y_1^{-1} &= \prod_{i=1}^k t_i^{h_i}.
\end{align*}
for some $k\ge0$, $t_i\in V$ and $g_i,h_i\in\bZ\setminus0$;
furthermore, $g_i$ and $h_i$ have the same sign for each $i$.
In particular, $x_1\cdot  y_1^{-1}$ is reduced
and so, the conclusion follows.
\ep

For each $u,v\in V$ and $x,y\in B(M)$, Lemma~\ref{l:uv} (1) implies the following  equivalences:
\begin{align*}
u^x = v^y &\Leftrightarrow u=v\text{ and }xy^{-1}\in Z(u) \\
 &\Leftrightarrow u=v\text{ and }\sigma(x)\sigma(y)^{-1}\in Z(u) \Leftrightarrow u^{\sigma x}=v^{\sigma y}.
 \end{align*}
\begin{align*}
\{u^x , v^y\}\in E(\Gamma^e)
 &\Leftrightarrow u\ne v\text{ and }\left[u^x,v^y\right]=1 \\
 &\Leftrightarrow  u\ne v,\ [u,v]=1,\ u^x=u^g,\ v^y = v^g\text{ for some }g\in A(\Gamma) \\
 &\Leftrightarrow  u\ne v,\ [u,v]=1,\ xy^{-1}\in Z(u)Z(v)\\
&\Leftrightarrow  u\ne v,\ [u,v]=1,\ \sigma(x)\sigma(y)^{-1}\in Z(u)Z(v)  \\
&\Leftrightarrow \{u^{\sigma x},v^{\sigma y}\}\in E(\Gamma^e).
 \end{align*}
So the map $\sigma$ satisfies the condition (i) of Lemma~\ref{l:sigma}, as desired.

\section*{Acknowledgements}
The authors thank Ilya Kapovich, Sang-jin Lee and Bertold Wiest for helpful discussions.

\bibliographystyle{amsplain}
\bibliography{ref}

\end{document}